\documentclass[reqno,12pt,a4paper]{amsart}

\usepackage{mathrsfs}
\usepackage{amssymb}
\usepackage{amsfonts}
\usepackage{latexsym}
\usepackage{amsthm}

\usepackage{graphicx}
\def\lb{\label}

\newcommand{\er}[1]{\textrm{(\ref{#1})}}

\begin{document}


\renewcommand{\theequation}{\arabic{section}.\arabic{equation}}
\theoremstyle{plain}
\newtheorem{theorem}{\bf Theorem}[section]
\newtheorem{lemma}[theorem]{\bf Lemma}
\newtheorem{corollary}[theorem]{\bf Corollary}
\newtheorem{proposition}[theorem]{\bf Proposition}
\newtheorem{definition}[theorem]{\bf Definition}
\newtheorem{remark}[theorem]{\it Remark}

\def\a{\alpha}  \def\cA{{\mathcal A}}     \def\bA{{\bf A}}  \def\mA{{\mathscr A}}
\def\b{\beta}   \def\cB{{\mathcal B}}     \def\bB{{\bf B}}  \def\mB{{\mathscr B}}
\def\g{\gamma}  \def\cC{{\mathcal C}}     \def\bC{{\bf C}}  \def\mC{{\mathscr C}}
\def\G{\Gamma}  \def\cD{{\mathcal D}}     \def\bD{{\bf D}}  \def\mD{{\mathscr D}}
\def\d{\delta}  \def\cE{{\mathcal E}}     \def\bE{{\bf E}}  \def\mE{{\mathscr E}}
\def\D{\Delta}  \def\cF{{\mathcal F}}     \def\bF{{\bf F}}  \def\mF{{\mathscr F}}
\def\c{\chi}    \def\cG{{\mathcal G}}     \def\bG{{\bf G}}  \def\mG{{\mathscr G}}
\def\z{\zeta}   \def\cH{{\mathcal H}}     \def\bH{{\bf H}}  \def\mH{{\mathscr H}}
\def\e{\eta}    \def\cI{{\mathcal I}}     \def\bI{{\bf I}}  \def\mI{{\mathscr I}}
\def\p{\psi}    \def\cJ{{\mathcal J}}     \def\bJ{{\bf J}}  \def\mJ{{\mathscr J}}
\def\vT{\Theta} \def\cK{{\mathcal K}}     \def\bK{{\bf K}}  \def\mK{{\mathscr K}}
\def\k{\kappa}  \def\cL{{\mathcal L}}     \def\bL{{\bf L}}  \def\mL{{\mathscr L}}
\def\l{\lambda} \def\cM{{\mathcal M}}     \def\bM{{\bf M}}  \def\mM{{\mathscr M}}
\def\L{\Lambda} \def\cN{{\mathcal N}}     \def\bN{{\bf N}}  \def\mN{{\mathscr N}}
\def\m{\mu}     \def\cO{{\mathcal O}}     \def\bO{{\bf O}}  \def\mO{{\mathscr O}}
\def\n{\nu}     \def\cP{{\mathcal P}}     \def\bP{{\bf P}}  \def\mP{{\mathscr P}}
\def\r{\rho}    \def\cQ{{\mathcal Q}}     \def\bQ{{\bf Q}}  \def\mQ{{\mathscr Q}}
\def\s{\sigma}  \def\cR{{\mathcal R}}     \def\bR{{\bf R}}  \def\mR{{\mathscr R}}
\def\S{\Sigma}  \def\cS{{\mathcal S}}     \def\bS{{\bf S}}  \def\mS{{\mathscr S}}
\def\t{\tau}    \def\cT{{\mathcal T}}     \def\bT{{\bf T}}  \def\mT{{\mathscr T}}
\def\f{\phi}    \def\cU{{\mathcal U}}     \def\bU{{\bf U}}  \def\mU{{\mathscr U}}
\def\F{\Phi}    \def\cV{{\mathcal V}}     \def\bV{{\bf V}}  \def\mV{{\mathscr V}}
\def\P{\Psi}    \def\cW{{\mathcal W}}     \def\bW{{\bf W}}  \def\mW{{\mathscr W}}
\def\o{\omega}  \def\cX{{\mathcal X}}     \def\bX{{\bf X}}  \def\mX{{\mathscr X}}
\def\x{\xi}     \def\cY{{\mathcal Y}}     \def\bY{{\bf Y}}  \def\mY{{\mathscr Y}}
\def\X{\Xi}     \def\cZ{{\mathcal Z}}     \def\bZ{{\bf Z}}  \def\mZ{{\mathscr Z}}
\def\O{\Omega}

\newcommand{\gA}{\mathfrak{A}}
\newcommand{\gB}{\mathfrak{B}}
\newcommand{\gC}{\mathfrak{C}}
\newcommand{\gD}{\mathfrak{D}}
\newcommand{\gE}{\mathfrak{E}}
\newcommand{\gF}{\mathfrak{F}}
\newcommand{\gG}{\mathfrak{G}}
\newcommand{\gH}{\mathfrak{H}}
\newcommand{\gI}{\mathfrak{I}}
\newcommand{\gJ}{\mathfrak{J}}
\newcommand{\gK}{\mathfrak{K}}
\newcommand{\gL}{\mathfrak{L}}
\newcommand{\gM}{\mathfrak{M}}
\newcommand{\gN}{\mathfrak{N}}
\newcommand{\gO}{\mathfrak{O}}
\newcommand{\gP}{\mathfrak{P}}
\newcommand{\gQ}{\mathfrak{Q}}
\newcommand{\gR}{\mathfrak{R}}
\newcommand{\gS}{\mathfrak{S}}
\newcommand{\gT}{\mathfrak{T}}
\newcommand{\gU}{\mathfrak{U}}
\newcommand{\gV}{\mathfrak{V}}
\newcommand{\gW}{\mathfrak{W}}
\newcommand{\gX}{\mathfrak{X}}
\newcommand{\gY}{\mathfrak{Y}}
\newcommand{\gZ}{\mathfrak{Z}}

\def\ve{\varepsilon}   \def\vt{\vartheta}    \def\vp{\varphi}    \def\vk{\varkappa}

\def\Z{{\mathbb Z}}    \def\R{{\mathbb R}}   \def\C{{\mathbb C}}    \def\K{{\mathbb K}}
\def\T{{\mathbb T}}    \def\N{{\mathbb N}}   \def\dD{{\mathbb D}}


\def\la{\leftarrow}              \def\ra{\rightarrow}            \def\Ra{\Rightarrow}
\def\ua{\uparrow}                \def\da{\downarrow}
\def\lra{\leftrightarrow}        \def\Lra{\Leftrightarrow}


\def\lt{\biggl}                  \def\rt{\biggr}
\def\ol{\overline}               \def\wt{\widetilde}
\def\no{\noindent}


\let\ge\geqslant                 \let\le\leqslant
\def\lan{\langle}                \def\ran{\rangle}
\def\/{\over}                    \def\iy{\infty}
\def\sm{\setminus}               \def\es{\emptyset}
\def\ss{\subset}                 \def\ts{\times}
\def\pa{\partial}                \def\os{\oplus}
\def\om{\ominus}                 \def\ev{\equiv}
\def\iint{\int\!\!\!\int}        \def\iintt{\mathop{\int\!\!\int\!\!\dots\!\!\int}\limits}
\def\el2{\ell^{\,2}}             \def\1{1\!\!1}
\def\sh{\sharp}
\def\wh{\widehat}
\def\bs{\backslash}

\def\where{\mathop{\mathrm{where}}\nolimits}
\def\all{\mathop{\mathrm{all}}\nolimits}
\def\as{\mathop{\mathrm{as}}\nolimits}
\def\Area{\mathop{\mathrm{Area}}\nolimits}
\def\arg{\mathop{\mathrm{arg}}\nolimits}
\def\const{\mathop{\mathrm{const}}\nolimits}
\def\det{\mathop{\mathrm{det}}\nolimits}
\def\diag{\mathop{\mathrm{diag}}\nolimits}
\def\diam{\mathop{\mathrm{diam}}\nolimits}
\def\dim{\mathop{\mathrm{dim}}\nolimits}
\def\dist{\mathop{\mathrm{dist}}\nolimits}
\def\Im{\mathop{\mathrm{Im}}\nolimits}
\def\Iso{\mathop{\mathrm{Iso}}\nolimits}
\def\Ker{\mathop{\mathrm{Ker}}\nolimits}
\def\Lip{\mathop{\mathrm{Lip}}\nolimits}
\def\rank{\mathop{\mathrm{rank}}\limits}
\def\Ran{\mathop{\mathrm{Ran}}\nolimits}
\def\Re{\mathop{\mathrm{Re}}\nolimits}
\def\Res{\mathop{\mathrm{Res}}\nolimits}
\def\res{\mathop{\mathrm{res}}\limits}
\def\sign{\mathop{\mathrm{sign}}\nolimits}
\def\span{\mathop{\mathrm{span}}\nolimits}
\def\supp{\mathop{\mathrm{supp}}\nolimits}
\def\Tr{\mathop{\mathrm{Tr}}\nolimits}
\def\BBox{\hspace{1mm}\vrule height6pt width5.5pt depth0pt \hspace{6pt}}


\newcommand\nh[2]{\widehat{#1}\vphantom{#1}^{(#2)}}
\def\dia{\diamond}

\def\Oplus{\bigoplus\nolimits}



\def\qqq{\qquad}
\def\qq{\quad}
\let\ge\geqslant
\let\le\leqslant
\let\geq\geqslant
\let\leq\leqslant
\newcommand{\ca}{\begin{cases}}
\newcommand{\ac}{\end{cases}}
\newcommand{\ma}{\begin{pmatrix}}
\newcommand{\am}{\end{pmatrix}}
\renewcommand{\[}{\begin{equation}}
\renewcommand{\]}{\end{equation}}
\def\eq{\begin{equation}}
\def\qe{\end{equation}}
\def\[{\begin{equation}}
\def\bu{\bullet}

\title[{Estimates for solutions of KDV in terms of action variables
}] {Estimates  for solutions of KDV on the phase space of
periodic distributions  in terms of action variables }

\date{\today}
\address{School of Math., Cardiff University.
Senghennydd Road, CF24 4AG Cardiff, Wales, UK.
email \ KorotyaevE@cf.ac.uk,
{\rm Partially supported by EPSRC grant EP/D054621.}}
\author[Evgeny L. Korotyaev]{Evgeny L. Korotyaev}


\begin{abstract}
We consider the KdV equation on the Sobolev space of periodic distributions. We obtain estimates of the solution of the KdV in terms of action variables.
\end{abstract}

\maketitle


\section{Introduction  and main results}
\setcounter{equation}{0}

Consider the KdV equation
\[
\pa_t\p=-\p_{xxx}'''+6\p\p_x'
\]
on the Sobolev space of zero-meanvalue  1-periodic distributions
$H_{-1}=\{\p=q': q\in H\}$, where  the real Hilbert space $H=H_0$ consists
 of zero-meanvalue
functions $q\in L^2(\T)$,  $\T=\R/\Z$.  The space $H_{-1}$ is  equipped
with the norm $\|\p\|_{{-1}}^2=\|q\|^2=\int_0^1 q^2(x)dx$ for $\p=q'\in H_{-1}$.
The initial value problem for KdV in the  phase space of
periodic distributions was solved by Kappeler and Topalov \cite{KT},
 see also \cite{B}, \cite{CT}.
That  problem  in various  Sobolev spaces
was studied by many authors, see references in  \cite{B}, \cite{CT}, \cite{KT}.
The action-angle variables for the periodic KdV are studied by Veselov-Novikov
 \cite{VN}, Kuksin \cite{Ku}, Kappeler-P\"oschel \cite{KP}. The action-angle
 variables for the case $\p\in H_{-1}$ were  constructed  by
 Kappeler-M\"ohr-Topalov \cite{KMT} and were  essentially used
 in \cite{KT}.
In  Sobolev
spaces  estimates  for the  potential $\psi$  and for  the KdV--Hamiltonian
  in terms of
the action variables  were obtained by Korotyaev \cite{K2}.

We describe the motivation of the present paper.
Introduce  the real Hilbert spaces
$\ell^2_{m}, m\in \R$,  of  sequences $(f_n)_1^{\iy }$, equipped
with the norm $\|f\|_m^2=\sum _{n\ge 1}(2\pi n)^{2m }f_n^2$. Recall that the
 KDV  equation on $H_{-1}$ admits action-angle variables
$A_n\ge 0, \f_n\in [0,2\pi), n\ge 1$  such that (see \cite{KT}):

\noindent
(1) for each $q\in H_{-1}$ there exist actions $A_n\ge 0$ such that $\sum_{n\ge 1}{A_n\/ n}<\iy$ and angles
 $\f_n\in [0,2\pi), n\ge 1$.

\noindent
(2) The mapping $\P: H_{-1}\to \ell_{-{1\/2}}^2\os \ell_{-{1\/2}}^2$
 given by $q\to \P(q)=(|A_n|^{1\/2}\cos \f_n,  |A_n|^{1\/2}\sin \f_n)_1^\iy$
 is a real analytic isomorphism between $H_{-1}$ and $\ell_{-{1\/2}}^2\os \ell_{-{1\/2}}^2$.

\noindent
(3) This mapping is symplectic.

For  the case of periodic distributions no  estimates for the potential $q$ in terms of action
variables $A_n$ were known.   Our main  goal in this paper is to obtain them.

We recall some results on the action variables  for the KdV.
 If $\p\in H_0$, then the following identity hold true
(see  \cite{MT}, \cite{K4}):
\[
\lb{ipq}
\|\p\|^2=4\sum_{n\ge 1} (\pi n) A_n.
\]
Moreover, if $\p'\in H$, then the Hamiltonian $\mH(\p)={1\/2}\int_0^1({\p'}^2+2\p^3)dx$
 obey the   following estimates:
\[
\lb{eHa}
8P_3-8P_1P_{-1}\le \mH(\p) \le 8P_3,\qqq
\]
where $P_j=\sum_{n\ge 1} (\pi n)^j A_n, j\in \R$, see \cite{K6}
and see \er{dA} for  definition of the actions  $A_n$.

We formulate our main result.

\begin{theorem}\lb{T1}
Let  $\p\in H_{-1}$ and let $P_{-1}=\sum_{n\ge 1}{A_n\/\pi n}$. Then the following
 estimates hold true:
\[
\lb{qA} \|\p\|_{-1}^2\le 3P_{-1}(1+P_{-1}),
\]
\[
\lb{Aq} P_{-1}\le \|\p\|_{-1}^2(1+\|\p\|_{-1})^{3\/2}.
\]
\end{theorem}

\no {\bf Conjecture.} The estimate \er{qA} is sharp.
It means that  an estimate $\|\p\|_{-1}^2\le
CP_{-1}(1+P_{-1})^\b$  with some $C>0$  and  $\b<1$ is not
correct.

We formulate a simple corollary, which follows directly
 from estimate \er{qA}, \er{Aq}.

\begin{corollary}
\lb{T2}
Let $\psi(x,t)$ be a solution of (1.1) such that
 $\p(\cdot,0)\in H_{-1}$. Then for all time $t$ the following estimates hold true:
\[
\lb{co}
\|\p(\cdot,t)\|_{-1}\le 3\|\p(\cdot,0)\|_{-1}(1+\|\p(\cdot,0)\|_{-1})^{5\/2},
\]
\[
\lb{coo}
\|\p(\cdot,0)\|_{-1}\le 14\max \{\|\p(\cdot,t)\|_{-1}, \|\p(\cdot,t)\|_{-1}^{5\/2}\}.
\]

\end{corollary}

{\bf Example. } We now discuss relation of estimates (1.6) with the inverse cascade of
energy in the KdV equation. Let an initial condition $\p(\cdot,0)\in H_0$ satisfies
\[
\lb{c1}
\|\p(\cdot,0)\|_{-1}=\ve\in [0,1/4],\qqq \qqq \|\p(\cdot,0)\|=C=\const.
\]
 Then for any $N\ge 1$ and every $t$  estimate \er{co} yields
\[
\lb{c2}
\|\p(\cdot,t)\|_{-1}\le 6\ve,\qqq \qqq \|\cP_N\p(\cdot,t)\|\le 6(2\pi N)\ve,
\]
where $\cP_N f, f\in H_0$ is given by
$
\cP_N f=\sum _{|n|\le N}e^{i2\pi nx}\int_0^1 f(s)e^{-i2\pi ns}ds.
$
Let in addition, $\d=6(2\pi N)\ve$ be small enough. Then
\er{c2} gives
\[
\lb{c3}
\|(I-\cP_N)\p(\cdot,t)\|^2\ge C^2-\d^2,\qqq
\|\cP_N\p(\cdot,t)\|\le \d, \qq any \qq t\ge 0.
\]
Thus we deduce that in our case the inverse cascade of energy is impossible.
It means that if the initial condition is
such that $\|\p(\cdot,0)\|=1$ and
$\|\cP_{N_0}\p(\cdot,0)\|=0$ for some $N_0\gg N$, then
$\|\cP_N\p(\cdot,t)\|\le {6N\/N_0}$ will be small for all time $t$, since
$\|\p(\cdot,0)\|_{-1}\le {1\/2\pi N}$. That is, if the energy of a solution
was initially concentrated in high modes, then a substantial part of
 the energy cannot flow to low modes.

Note that the function $\psi(x,0)$  with the  property \er{c1} may be a finite
trigonometric polynomial.

\section{Proof of the main theorem}
\setcounter{equation}{0}

Our main ingredients  to study the KdV equation (similar to \cite{KT})
 are the spectral properties of  the Schr\"odinger operator $T=-{d^2\/dx^2}+\p+q_0$, where $\p\in H_{-1}$ is
a 1-periodic distribution with zero mean-value  and $q_0\in \R$ is a constant.
 It is well known \cite{K3}  that the spectrum of $T$  is
absolutely continuous and consists of intervals
$\gS_n=[\l^+_{n-1},\l^-_n ]$,  where $\l^+_{n-1}<\l^-_n \le  \l^+_{n},\
n\ge 1$. We take a constant $q_0$ such that $\l_0^+=0$.  The intervals $\gS_{n}$ and
$\gS_{n+1}$  are separated by the gap $\g_n=(\l^-_n,\l^+_n )$.  If a gap degenerates, that is $\g_n=\es$,
then the corresponding segments $\gS_{n} $ and $\gS_{n+1}$ merge.
 The sequence $\l_0^+<\l_1^-\le \l_1^+\ <\dots$
 is the spectrum of the equation $-y''+(\p+q_0)y=\l y$ with the 2-periodic
 boundary conditions, i.e.,  $y(x+2)=y(x), x\in \R$.
If $\l_n^-=\l_n^+$ for some $n$, then this number $\l_n^{\pm}$
 is the double eigenvalue of this equation with the 2-periodic
 boundary conditions. The lowest  eigenvalue $\l_0^+$ is always simple and the
corresponding eigenfunction is 1-periodic. The eigenfunctions,
corresponding to the eigenvalue $\l_n^{\pm}$, are 1-periodic,
when $n$ is even and  are antiperiodic,  i.e., $y(x+1)=-y(x),\ \
x\in\R$, when $n$ is odd.

We can not introduce the standard fundamental solutions for the
operator $T$ since the perturbation $\p\in H_{-1}$ is very strong. But
 we can do this using another representation of $T$ given by $T=\mU T_w\mU^{-1}$.
Here $T_w$ is the self-adjoint periodic operator acting in $L^2(\R,w^2(x)dx)$ and given by
\[
\lb{Tw} \qqq \qqq T_wf=-{1\/w^2}(w^2f')'=-f''-2pf',  \
 w(x)=e^{\int_0^x p(s)ds}, \ \ p\in H.
\]
$\mU$ is the unitary transformation $\mU:L^2(\R,w^2 dx)\to L^2(\R,dx)$, given by the
 multiplication by $w$.
Note that
\[
\lb{1}
T=-{d^2\/dx^2}+q'+q_0\ge 0,\qqq q'=p'(x)+p^2(x)-\|p\|^2,\qqq
 q_0=\|p\|^2=\int_0^1 p^2(x)dx,
\]
where $\p=q'$ is  a 1-periodic
potential (distribution).  Thus, if $p'\in H$, then $T_w$
corresponds to the Hill  operator $T$ with $L^2$-potential.
The operator $T_w$ is well studied, see \cite{K5}
 and references therein. In fact the direct spectral problem for $T_w$
is equivalent to that for $T$ \cite{K3}.

The operator $T_w$ has the standard fundamental solutions
$\vp(x,\l), \vt(x,\l)$, which satisfy the  equation $-y''-2py'=\l y,
\ \l\in \C$ and the conditions $\vp'(0,\l)=\vt(0,\l)=1,\
\vp(0,\l)=\vt'(0,\l)=0$. Here and below we use the notation $f'={\pa
\/\pa x}f$. Introduce the Lyapunov function
$\D(\l)={1\/2}(\vp'(1,\l)+\vt(1,\l))$.  Note
 that $\D(\l_{n}^{\pm})=(-1)^n,\ n\ge 1,$ and that for
each $n\ge 1$ there exists a unique point  $\l_n\in [\l^-_n,\l^+_n]$ such that  $\D'(\l_n)=0$.

Now we recall results, crucial for the present paper.
For each $\p\in H_{-1}$ there exists a unique conformal mapping (the
quasimomentum) $k:\cZ\to \cK(h) $ with asymptotics $k(z)=z+o(1)$ as
$|z|\to \iy$ (see Fig. 1 and 2) and such that (see \cite{K3})
\begin{multline}
\lb{dequ} \cos k(z)= \D(z),\ \  z\in \cZ =\C\sm\cup \ol g_n,\qq
g_n=(z^-_n,z^+_n )=-g_{-n},\qq z_n^{\pm}=\sqrt{\l_n^{\pm}}\ge 0,\qq n\ge 1,
\\
\cK(h)=\C\sm\cup \G_n ,\ \  \  \G_n=(\pi n-ih_n,\pi n+ih_n),\qqq
h_0=0,\qqq h_n=h_{-n}\ge 0,
\qqq
\\
h_n\ge 0 \qqq {\rm is \ defined \ by  \ the \ equation}\qqq  \cosh
h_n = (-1)^n\D(\l_n) \ge 1.
\end{multline}
Here $g_0=\es$ and $\G_n$ is the vertical cut, $z_n=\sqrt{\l_n}\in
[z_n^-,z_n^+], n\ge 1$,  $\D'(z_n^2)=0$. Moreover, we have $(h_n)_1^\iy\in \ell^2$ iff $\p\in H_{-1}$ (and $(nh_n)_1^\iy\in
\ell^2$ iff $\p\in H$), see \cite{K3}, \cite{K1}.

Due to  \cite{MO1}, the quantities  $v=\Im k(z)$ and $u=\Re k(z), z\in \cZ$, possess the
following properties:

\no 1) {\it $v(z)\ge \Im z>0$ and $v(z)=-v(\ol z)$ for all $z\in
\C_+=\{\Im z>0\}$.

\no 2) $v(z)=0$ for all $z\in \s_n=[z^+_{n-1},z^-_n]=-\s_{-n}, n\ge 1$.

\no 3)  If some $g_n\ne \es, n\in \Z$, then the function $v(z+i0)>0$
for all $z\in g_n$, and $v(z+i0)$ has a maximum at $z_n\in g_n$ such
that $\D'(z_n^2)=0$ and $v(z_n+i0)=h_n>0, v'(z_n)=0$,
 and
\[
\lb{prq} v(z+i0)=-v(z-i0)>v_n(z)=|(z-z_n^-)(z-z_n^+)|^{1\/2}>0, \qqq v''(z+i0)<0,
\]
for all $ z\in
g_n\ne \es$, see Fig. 3.

\no 4) $u'(z)>0$ on  $\R\sm \cup \ol g_n$ and  $u(z)=\pi n$ for
all $z\in g_n\ne \es, n\in \Z$.

\no 5) The function $k(z)$ maps a horizontal cut (a "gap" )
$\ol g_n$  onto a vertical cut $\G_n$  and a spectral band
$\s_n$ onto the segment $[\pi (n-1), \pi n]$ for all $\pm n\in\N$. }

The heights $h_n,\  n\ge 1$ are so-called Marchenko-Ostrovski
parameters \cite{MO1}.
In spirit, such result goes back to the classical Hilbert Theorem
(for a finite number of cuts, see e.g. \cite{J}) in the conformal
mapping theory. A similar theorem for the Hill operator is
technically more complicated (there is a infinite number of cuts)
and was proved by Marchenko-Ostrovski \cite{MO1} for the case $\p \in H$.  For additional properties of the conformal mapping we also refer to our previous papers \cite{K1}-\cite{K6}.
Note that the inverse problems for the operator $H$ with $\p\in H_{-1}$ in terms of the Marchenko-Ostrovski parameters
$h_n, n\ge 1$  and gap-lengths were solved by Korootyev in \cite{K3}.

For the sake of the reader, we briefly recall the results existing in the
literature about estimates.
In the case $h=(h_n)_1^\iy$ and $\p\in H$  Marchenko and Ostrovki [MO1-2] obtained  the estimates:
$
\|\p\|\le C(1+\sup_{n\ge 1}h_n)\|h\|_1 $ and
$\|h\|_1\le C\|\p\|\exp (\ C_1\|\p\|\ )$
 for some absolute constants $C, C_1$. These estimates are not sharp since they used the Bernstein inequality.
Using the harmonic measure argument  Garnett and Trubowitz \cite{GT}
obtained $\|\g\|\le (4+\|h\|_1)\|h\|_1$ for the case $\p\in H$ and $\g=(|\g_n|)_1^\iy$, where $|\g_n|$ is a gap length.
Using the conformal mapping theory,
Korotyaev \cite{K1}-\cite{K6} obtained estimates  of potentials (and the Hamiltonian of the KDV) in terms of  gap lengths, actions variables, effective masses, the heights $h=(h_n)_1^\iy$ for large class of potentials. In fact in order to get new estimates  new results from the comformal mapping theory were obtained. Note that estimates simplify the proof for the inverse
spectral theory, see \cite{KK}, \cite{K3}. We recall only few results from these estimates:

\no {\it I). Let $\p\in H_{-1}$. Then the following estimates hold true (see \cite{K3}):
\[
\lb{TA1} \|\g\|_{-1}\le \sqrt 2 \|\p\|_{-1}(1+\|\p\|_{-1}),\ \ \ \ \ \|\p\|_{-1}\le 8\pi \|\g\|_{-1}(1+\|\g\|_{-1}),
\]
\[
\lb{TA4} {\sqrt\pi\/ \sqrt 8}\|\p\|_{-1}\le \|h\|_0\le {\pi\/
2}\|\p\|_{-1}(1+\|\p\|_{-1})^{1\/2}.
\]
\no II) If $\p\in H$, then the following estimates hold true
(see \cite{K1})}:
$$
\|\p\|\le 2\|\g\|_0(1+\|\g\|_0^{1\/3}),\qqq
\|\g\|_0\le 2\|\p\|(1+\|\p\|^{1\/3}).
$$

If $\p\in H_{-1}$, then the quasimomentum $k(\cdot)$ has asymptotics
$$
k(z)=z-{Q_0+o(1)\/z}\qqq as \qq z\to +i\iy,
$$
where $Q_0={1\/\pi}\int_\R v(z+i0)dz\ge 0$ and  $p$ (defined in \er{1})  satisfy the identities from \cite{K5}:
\[
\lb{HQ}
Q_0={1\/\pi}\int_\R v(z+i0)dz={1\/2\pi}\int\!\!\int_\C|z'(k)-1|^2dudv={\|p\|^2\/2}, \qq k=u+iv.
\]

Due to \cite{FM} we define the action $A_n, n\ge 1$ by
\[
\lb{dA}
A_n={(-1)^{n+1}2\/\pi}\int_{\g_n}{\l\D'(\l)d\l\/|\D^2(\l)-1|^{1\/2}}\ge 0.
\]
We rewrite $A_n$ in the more convenient form. The differentiation of
$\D(z^2)=\cos k(z)$ gives $k'(z)=-{\D'(z^2)2z\/\sin k(z)}$, which together with $\sin
k(z)=\sqrt{1-\D^2(z^2)}$ yield
$$
A_n=-{1\/i\pi}\int_{c_n}z^2{\D'(z^2)2z\/\sin
k(z)}dz={1\/i\pi}\int_{c_n}z^2k'(z)dz=-
{2\/i\pi}\int_{c_n}zk(z)dz={4\/\pi}\int_{g_n}zv(z+i0)dz\ge 0,
$$
which gives
\[
\lb{iA}
A_n={4\/\pi}\int_{g_n}zv(z+i0)dz\ge 0.
\]

\begin{figure}
\tiny \unitlength=1mm \special{em:linewidth 0.4pt}
\linethickness{0.4pt}
\begin{picture}(120.67,34.33)
\put(20.33,21.33){\line(1,0){100.33}}
\put(70.33,10.00){\line(0,1){24.33}}
\put(69.00,19.00){\makebox(0,0)[cc]{$0$}}
\put(120.33,19.00){\makebox(0,0)[cc]{$\Re z$}}
\put(67.00,33.67){\makebox(0,0)[cc]{$\Im z$}}
\put(81.33,21.33){\linethickness{2.0pt}\line(1,0){9.67}}
\put(100.33,21.33){\linethickness{2.0pt}\line(1,0){4.67}}
\put(116.67,21.33){\linethickness{2.0pt}\line(1,0){2.67}}
\put(60.00,21.33){\linethickness{2.0pt}\line(-1,0){9.33}}
\put(40.00,21.33){\linethickness{2.0pt}\line(-1,0){4.67}}
\put(24.33,21.33){\linethickness{2.0pt}\line(-1,0){2.33}}
\put(81.67,24.00){\makebox(0,0)[cc]{$z_1^-$}}
\put(91.00,24.00){\makebox(0,0)[cc]{$z_1^+$}}
\put(100.33,24.00){\makebox(0,0)[cc]{$z_2^-$}}
\put(105.00,24.00){\makebox(0,0)[cc]{$z_2^+$}}
\put(115.33,24.00){\makebox(0,0)[cc]{$z_3^-$}}
\put(120.00,24.00){\makebox(0,0)[cc]{$z_3^+$}}
\put(59.33,24.00){\makebox(0,0)[cc]{$-z_1^-$}}
\put(50.67,24.00){\makebox(0,0)[cc]{$-z_1^+$}}
\put(40.33,24.00){\makebox(0,0)[cc]{$-z_2^-$}}
\put(34.67,24.00){\makebox(0,0)[cc]{$-z_2^+$}}
\put(26.00,24.00){\makebox(0,0)[cc]{$-z_3^-$}}
\put(19.50,24.00){\makebox(0,0)[cc]{$-z_3^+$}}
\end{picture}
\caption{Domain $\cZ=\C\sm\cup g_n$, where
 $z=\sqrt{\l}$ and momentum gaps $g_n=(z_n^-,z_n^+)$}
\lb{z}
\end{figure}
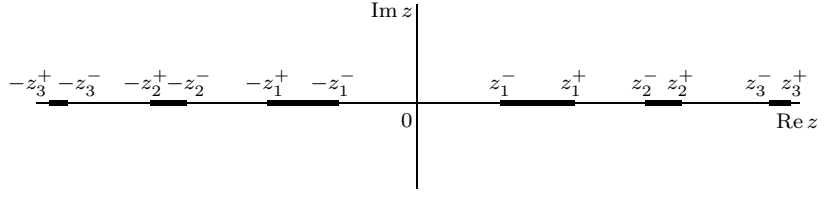
%
%
%
\begin{figure}
\tiny \unitlength=1mm \special{em:linewidth 0.4pt}
\linethickness{0.4pt}
\begin{picture}(120.67,34.33)
\put(20.33,20.00){\line(1,0){102.33}}
\put(71.00,7.00){\line(0,1){27.00}}
\put(70.00,18.67){\makebox(0,0)[cc]{$0$}}
\put(124.00,18.00){\makebox(0,0)[cc]{$\Re k$}}
\put(67.00,33.67){\makebox(0,0)[cc]{$\Im k$}}
\put(87.00,15.00){\linethickness{2.0pt}\line(0,1){10.}}
\put(103.00,17.00){\linethickness{2.0pt}\line(0,1){6.}}
\put(119.00,18.00){\linethickness{2.0pt}\line(0,1){4.}}
\put(56.00,15.00){\linethickness{2.0pt}\line(0,1){10.}}
\put(39.00,17.00){\linethickness{2.0pt}\line(0,1){6.}}
\put(23.00,18.00){\linethickness{2.0pt}\line(0,1){4.}}
\put(85.50,18.50){\makebox(0,0)[cc]{$\pi$}}
\put(54.00,18.50){\makebox(0,0)[cc]{$-\pi$}}
\put(101.00,18.50){\makebox(0,0)[cc]{$2\pi$}}
\put(36.00,18.50){\makebox(0,0)[cc]{$-2\pi$}}
\put(117.00,18.50){\makebox(0,0)[cc]{$3\pi$}}
\put(20.00,18.50){\makebox(0,0)[cc]{$-3\pi$}}
\put(87.00,26.00){\makebox(0,0)[cc]{$\pi+ih_1$}}
\put(56.00,26.00){\makebox(0,0)[cc]{$-\pi+ih_1$}}
\put(103.00,24.00){\makebox(0,0)[cc]{$2\pi+ih_2$}}
\put(39.00,24.00){\makebox(0,0)[cc]{$-2\pi+ih_2$}}
\put(119.00,23.00){\makebox(0,0)[cc]{$3\pi+ih_3$}}
\put(23.00,23.00){\makebox(0,0)[cc]{$-3\pi+ih_3$}}
\end{picture}
\caption{$k$-plane and cuts $\G_n=(\pi n-ih_n,\pi n+ih_n), n\in\Z$}
\lb{k}
\end{figure}
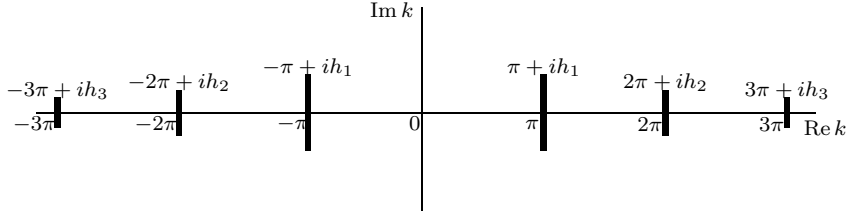
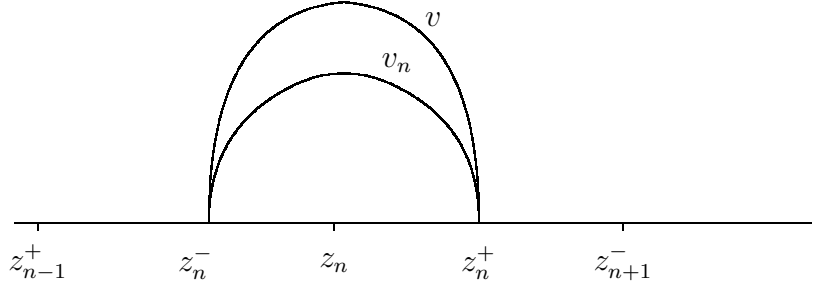
\begin{figure}
\unitlength 1mm 
\linethickness{0.4pt}
\ifx\plotpoint\undefined\newsavebox{\plotpoint}\fi 
\begin{picture}(119.75,66.5)(0,0)
\put(15,17.25){\line(1,0){104.75}}
\qbezier(40.5,17.25)(40.5,30.375)(52.25,35.75)
\qbezier(52.25,35.75)(58.25,38.375)(64.25,35.75)
\qbezier(64.25,35.75)(76.17,30.375)(76,17.25)
\qbezier(40.5,17.25)(40.5,45.125)(58.25,46.5)
\qbezier(76,17.25)(76.17,45.125)(58.25,46.5)
\put(38.75,12.00){\makebox(0,0)[cc]{$z_n^-$}}
\put(76,12){\makebox(0,0)[cc]{$z_n^+$}}
\put(65.25,38.5){\makebox(0,0)[cc]{$v_n$}}
\put(70,44.25){\makebox(0,0)[cc]{$v$}}
\put(18,17.25){\line(0,-1){1.00}} \put(95,17.25){\line(0,-1){1.00}}
\put(57,17.25){\line(0,-1){1.00}}
\put(18,12.00){\makebox(0,0)[cc]{$z_{n-1}^+$}}
\put(95,12.00){\makebox(0,0)[cc]{$z_{n+1}^-$}}
\put(57,12.00){\makebox(0,0)[cc]{$z_{n}$}}
\end{picture}
\caption{The graph of $v(z+i0), \ z\in g_n\cup \s_n\cup \s_{n+1}$
and $|h_n|=v(z_n+i0)>0$} \lb{grafv}
\end{figure}

Below we need
results about the Riccati mapping (see Theorem 1.2 in \cite{K3}).

\begin{theorem}   \lb{Tri}
The Riccati map $R: H\to H$ given by $p\to q=R(p), q'=p'(x)+p^2(x)-\|p\|^2$
is a real analytic isomorphism of $H$ onto itself. Moreover, the following
estimates hold true:
\[
\lb{Tri1}
\|q\|\le \|p\|(1+2\|p\|),\ \ \ \ \
\]
\[
\lb{Tri2}
 \|p\| \le \sqrt 2\|q\|(1+2\|q\|). \ \ \  \ \
\]
\end{theorem}

In order to show \er{qA}, we need

\begin{lemma}   \lb{TRs}
 The following estimate holds true:
\[
\lb{TRs2} \|p\|^2\le \sum_{n\ge 1}{A_n\/\pi n}=P_{-1}.
\]

\end{lemma}
\no {\bf Proof.}  Using the following identity for $Q_0={1\/\pi}\int_\R v(z+i0)dz$
(see Theorem 2.3 from \cite{K4})
\[
Q_0={2\/\pi}\int_0^\iy {u(z)v(z)\/z}dz,
\]
 we obtain
$$
Q_0^2\le {2\/\pi}\int_{g_+}{uvdz\/z} {2\/\pi}\int_{g_+}{zvdz\/u}=
Q_0{2\/\pi}\int_{g_+}{zvdz\/u},\qqq g_+=\cup_{n\ge 1}g_n,
$$
which together with the identity for $A_n$
\er{iA} yields
$$
Q_0\le{2\/\pi}\int_{g_+}{zvdz\/u}={1\/2}\sum_{n\ge 1}{A_n\/\pi n},
$$
since $u|_{g_n}=\pi n$ and the identity \er{HQ} gives \er{TRs2}.
\BBox

We show the estimate \er{qA}. Using \er{Tri1}, \er{TRs2} we obtain
$$
\|\p\|_{-1}^2=\|q\|^2\le \|p\|^2(1+2\|p\|)^2\le \|p\|^2(1+\|p\|^2)
\le 5P_{-1}(1+P_{-1}),
$$
which gives \er{qA}.

We show the estimate \er{Aq}. The estimate $v|_{g_n}\le h_n$ and the identity for $A_n$
\er{iA} gives
$$
A_n={4\/\pi}\int_{g_n}zv(z)dz\le
{4h_n\/\pi}\int_{g_n}zdz={4h_n\/\pi}|\g_n|,
$$
and then
\[
\lb{pe} P_{-1}=\sum_{n\ge 1}{A_n\/\pi n}\le \sum_{n\ge
1}{4h_n\/\pi}{|\g_n|\/\pi n} \le {4\/\pi}\|h\|_0 \|\g\|_{-1}.
\]
 Substituting estimates \er{TA1}, \er{TA4}
into \er{pe} we obtain \er{Aq}. \BBox

\no {\bf Acknowledgments.}
\small
The author is grateful to Sergei Kuksin     (Ecole Polytechnique, Paris)  for stimulating discussions and useful comments.

\end{document}